\documentclass[12pt,psfig,reqno]{amsart}
\usepackage{amsfonts}
\usepackage{amscd}
\pdfoutput=1
\usepackage{amssymb,amsfonts,latexsym}
\usepackage{graphics,verbatim}
\usepackage{graphicx}
\usepackage[usenames]{color} 
\makeatletter

\renewcommand\section{\@startsection {section}{1}{\z@}%
                                   {-3.5ex \@plus -1ex \@minus -.2ex}%
                                   {2.3ex \@plus.2ex}%
                                   {\centering\normalfont\bf}}

 \numberwithin{equation}{section}
\numberwithin{equation}{section}

\numberwithin{equation}{section}

\theoremstyle{plain}

\begin{document}
\title{   Several Integral Estimates  and Some  Applications }
\author{  Xuejun  Zhang $^{*}$, Hongxin Chen, Min Zhou, Yuting Guo, Pengcheng Tang  }
\address{College of Mathematics and Statistics, Hunan Normal University, Changsha, Hunan 410006, China}

\email{xuejunttt@263.net}

\date{}
\thanks{$^*$ Corresponding author.\\
 The research is supported by the National
Natural Science Foundation of China (No. 11942109) and the Natural Science Foundation of Hunan Province(No. 2022JJ30369).}

\begin{abstract}  In this paper, the authors first consider the bidirectional estimates of several typical integrals. As some applications of these integral estimates,
 the authors investigate the pointwise multipliers  from the normal weight  general function space $F(p,\mu,s)$ to the normal weight Bloch type space $\mathcal{B_{\nu}}(B_{n})$ on the unit ball $B_{n}$ of $\mathbb{C}^{n}$, where $\mu$ and $\nu$ are two normal functions on $[0,1)$. For the special normal function $\displaystyle{\mu(r)=(1-r^{2})^{\alpha}\log^{\beta}\frac{e}{1-r^{2}}}$ \ ($\alpha>0$, \ $-\infty<\beta<\infty$), the authors give the necessary and sufficient conditions of pointwise multipliers from  $F(p,\mu,s)$ to $\mathcal{B_{\nu}}(B_{n})$ for all cases.

\end{abstract}
\maketitle
\section{\bf Introduction  \label{sect.1}}

\ \ \ In this paper, we call  $``P\asymp Q"$ if there exist two positive constants $A_{1}$ and $A_{2}$ such that $A_{1}P \leq Q\leq A_{2}P$. We  call  $``Q\lesssim P"$ (or $``Q \gtrsim P"$) if there exists a positive constant $A$ such that $ Q \leq AP$ (or $ Q \geq AP$).
\vskip2mm
Let $\mathbb{C}$ denote the set of complex numbers. Fixed a positive integer $n$, let $\mathbb{C}^{n}=\mathbb{C}\times\cdots\times\mathbb{C}$ denote the Euclidean space of complex dimension $n$. For $z=(z_{1},\cdots,z_{n})$ and $w=(w_{1},\cdots,w_{n})$ in  $\mathbb{C}^{n}$, the inner product of $z$ and $w$ is defined by $\langle z,w\rangle=z_{1}\overline{w_{1}}+\cdots+z_{n}\overline{w_{n}}$. The open unit ball in $\mathbb{C}^{n}$ is the set
$B_{n}=\{z\in \mathbb{C}^{n}: |z|=\sqrt{\langle z,z\rangle}<1\}$ \ (we write as $D$ when $n=1$), and the boundary of $B_{n}$ is the set $S_{n}=\{z\in \mathbb{C}^{n}: |z|=1\}$.
The class of
holomorphic functions and the space of bounded holomorphic functions on $B_{n}$  are denoted by $H(B_{n})$ and $H^{\infty}(B_{n})$,
 respectively. For $a\in B_{n}$, let $\varphi_{a}$ be the holomorphic automorphism of $B_{n}$ satisfying with $\varphi_{a}(0)=a$,  $\varphi_{a}(a)=0$ and $\varphi_{a}^{-1}=\varphi_{a}$. It follows from Lemma 1.3 in [1] that
$$
1-\langle \varphi_{a}(z),\varphi_{a}(w)\rangle=\frac{(1-|a|^{2})(1-\langle z,w\rangle)}{(1-\langle z,a\rangle)(1-\langle a,w\rangle)} \ \ \ (z,w\in B_{n}).\eqno{(1.1)}
$$
In particular, if $w=z $ or $w=0$, then there are
$$
1-|\varphi_{a}(z)|^{2}=\frac{(1-|a|^{2})(1-|z|^{2})}{|1-\langle z,a\rangle|^{2}}, \ \ \ 1-\langle \varphi_{a}(z), a\rangle=\frac{1-|a|^{2}}{1-\langle z,a\rangle}.\eqno{(1.2)}
$$
 Given  $r>0$,  the Bergman ball with $a$ as the center and $r$ as the radius is the set
$$
D(a,r)=\{z\in B_{n}: \  \beta(z,a)<r\}, \ \ \mbox{where} \ \ \beta(z,a)=\frac{1}{2}\log\frac{1+|\varphi_{a}(z)|}{1-|\varphi_{a}(z)|}.
$$

1996,  R. H. Zhao first introduced the general function space $F(p,q,s)$ on $D$ (see [2]).
 Soon after, the space $F(p,q,s)$ was extended to $B_{n}$. If $f\in F(p,q,s)$, then
 \begin{eqnarray*}
 &\;&
\sup_{w\in B_{n}}\int_{B_{n}}|\nabla f(z)|^{p}(1-|z|^{2})^{q}\log^{s}\frac{1}{|\varphi_{w}(z)|}\ dv(z)\\
&\;&\asymp \sup_{w\in B_{n}}(1-|w|^{2})^{s}\int_{B_{n}}\frac{(1-|z|^{2})^{q+s}|\nabla f(z)|^{p}}{|1-\langle z,w\rangle|^{2s}}\ dv(z),
 \end{eqnarray*}
 where $\displaystyle{\nabla f(z)=\left(\frac{\partial f}{\partial z_{1}}(z),\cdots, \frac{\partial f}{\partial z_{n}}(z)\right)}$ (see [3]).
 \vskip2mm
In recent years, there  have been a lot of  results related to $F(p,q,s)$ space, such as [2]-[18] etc. In order to consider the general function space in a broader field of vision, we first give the definition of normal function.

\vskip1mm {\bf Definition 1.1} \ A positive continuous function on $[0,1)$ is called normal if there exist constants $0< a\leq
b<\infty$ and $0\leq r_{0}<1$ such that
\vskip2mm
(i) \ $\displaystyle{\frac{\mu(r)}{(1-r^{2})^a}}$  is decreasing on $[r_{0},1)$; \ (ii) \ $\displaystyle{\frac{\mu(r)}{(1-r^{2})^b}}$ is increasing on $[r_{0},1)$.

\vskip2mm
For example,
\ $\displaystyle{\mu(r)=(1-r^{2})^{\alpha}\log^{\beta}\frac{e}{1-r^{2}}\left\{\log\log\frac{e^{2}}{1-r^{2}}\right\}^{\gamma}}$ \
($\alpha>0$, \ $-\infty<\beta, \ \gamma<\infty$)
 is this kind of  normal function.
 \vskip2mm
 Without losing generality,  let
$r_{0}=0$ in this paper.

\vskip2mm
Recently, S. L. Li ([4]) extended  $F(p,q,s)$  to a kind of abstract form as follows:

\vskip2mm
{\bf  Definition 1.2} \ Let $\mu$ be a normal function on $[0,1)$.  For $p>0$ and $s\geq 0$, if $f\in H(B_{n})$
$$
\mbox{and} \ \ ||f||_{p,\mu,s}=|f(0)|+\left\{\sup_{w\in B_{n}}\int_{B_{n}}\frac{(1-|w|^{2})^{s}|\nabla f(z)|^{p}}{|1-\langle z,w\rangle|^{2s}}\frac{\mu^{p}(|z|)}{1-|z|^{2}}\ dv(z)\right\}^{\frac{1}{p}}<\infty,
$$
then $f$ is said to belong to the normal weight general function space $F(p,\mu,s)$.
\vskip2mm
 It is easy to prove that $F(p,\mu,s)$ is a Banach space with the norm $||.||_{p,\mu,s}$ when $p\geq 1$, and $F(p,\mu,s)$ is a complete distance space with  $\rho(f,g)=||f-g||_{p,\mu,s}^{p}$ when $0<p<1$.
\vskip2mm
This paper also involves the normal weight Bloch type space and the weighted Bergman space. We give their definitions, respectively.
\vskip2mm
{\bf Definition 1.3} \ Let $\mu$ be a normal function on $[0,1)$. The  normal weight Bloch type space $\mathcal{B_{\mu}}(B_{n})$ consists of holomorphic function $f$ in  $ B_{n}$ such that
$$
||f||_{\mathcal{B_{\mu}}}=|f(0)|+\sup_{z\in B_{n}}\mu(|z|)|\nabla f(z)|<\infty.
$$

The space $\mathcal{B_{\mu}}(B_{n})$ is a Banach space under the norm $||.||_{\mathcal{B_{\mu}}}$. In particular, $\mathcal{B_{\mu}}(B_{n})$ is just the  $\alpha$-Bloch space  $\mathcal{B^{\alpha}}(B_{n})$  when $\mu(r)=(1-r^{2})^{\alpha} \ (\alpha>0)$.

\vskip2mm
{\bf Definition 1.4} \ For $\alpha>-1$ and \ $p>0$, the weighted Bergman space $A^{p}_{\alpha}(B_{n})$ consists of holomorphic function $f$ in  $ B_{n}$ such that \ $
\int_{B_{n}}|f(z)|^{p}\ dv_{\alpha}(z)<\infty$,
 where  $dv_{\alpha}(z)=c_{\alpha}(1-|z|^{2})^{\alpha}\ dv(z)$, \ $c_{\alpha}=\Gamma(n+1+\alpha)/n!\Gamma(\alpha+1)$ and  $v(B_{n})=1$.

\vskip2mm
It is known that  the application of integral estimate is a basic tool for the study of function space theory and operator theory. Next, we introduce several integral estimates.
\vskip2mm
 For a point in $B_{n}$,  W. Rudin gave the following Proposition in [19]:
\vskip2mm
{\bf Proposition A} \ Let $t>-1$ and $c$ be  real. Then the integrals
$$
I(z)=\int_{S_{n}}\frac{d\sigma(\xi)}{|1-\langle
\xi,z,\rangle|^{n+c}} \ \ \mbox{and} \ \ J(z)=\int_{ B_{n}}\frac{(1-|w|^{2})^{t}}{|1-\langle
z,w\rangle|^{n+1+t+c}}\ dv(w) \ \ (z\in B_{n})
$$
have the following asymptotic properties:
\vskip2mm
(1) \ $I(z)\asymp J(z)\asymp 1$ \ when $c<0$. \ (2) \ $\displaystyle{I(z)\asymp J(z)\asymp
\log\frac{e}{1-|z|^{2}}}$ \ when $c=0$.\vskip2mm (3) \ $\displaystyle{I(z)\asymp J(z)\asymp
\frac{1}{(1-|z|^{2})^{c}}}$ \ when $c>0$.
\vskip2mm
For two  points in $B_{n}$, the following Proposition B comes from Lemma 2.2 in [10] (case $k=0$) and Theorem 3.1 in [12] (case $k>0$), and the following Proposition C comes from Proposition 3.1 in [36].

\vskip2mm
{\bf Proposition B} \ Let \ $\delta>-1$, \ $r\geq 0$, \ $t\geq 0$ \ and \ $k\geq 0$.  Then  the integral
\begin{eqnarray*}
 &\;&
 J_{w,a}=\int_{B}\frac{(1-|z|^{2})^{\delta}}{|1-\langle z, w\rangle|^{t}\ |1-\langle z,a\rangle|^{ r}}\log^{k}\frac{e}{1-|z|^{2}} \ dv(z) \ \ (w,a\in B_{n})
 \end{eqnarray*}
has the following bidirectional estimates:
\vskip2mm (1) \ \ $ J_{w,a}\asymp 1 $ \ when \ $t+r-\delta<n+1$.
\vskip2mm (2) \ \ $ \displaystyle{J_{w,a}\asymp \log^{k+1}\frac{e}{|1-\langle
w,a\rangle|} }$ \ when \ $t+r-\delta=n+1$, \ $t>0$ \ and \ $r>0$.
\vskip2mm (3) \ \ $ \displaystyle{J_{w,a}\asymp\frac{1}{|1-\langle
w,a\rangle|^{r}}\log^{k}\frac{e}{1-|w|^{2}}\log\frac{e}{|1-\langle w,\varphi_{w}(a)\rangle|} }
$ \vskip2mm when \ $t-\delta=n+1>r-\delta$ \ and \ $r>0$.
\vskip2mm (4) \ \ $ \displaystyle{J_{w,a}\asymp\frac{1}{|1-\langle
w,a\rangle|^{t+r-\delta-n-1}} }\log^{k}\frac{e}{|1-\langle w,a\rangle|}
$ \vskip2mm when \ $t+r-\delta>n+1>\max\{r-\delta, \ t-\delta\}$.
\vskip2mm
(5) \ \ $ \displaystyle{J_{w,a}\asymp \frac{1}{|1-\langle
w,a\rangle|^{\delta+n+1}}\log^{k}\frac{e}{1-|w|^{2}}\log\frac{e}{|1-\langle w,\varphi_{w}(a)\rangle|}}$
\vskip2mm + \ $\displaystyle{\frac{1}{|1-\langle
w,a\rangle|^{\delta+n+1}}\log^{k}\frac{e}{1-|a|^{2}}\log\frac{e}{|1-\langle a,\varphi_{a}(w)\rangle|}} $ \ when \ $t-\delta=n+1=r-\delta$.
\vskip2mm
(6) \ \
$ \displaystyle{I_{w,a}\asymp \frac{1}{(1-|w|^{2})^{ t-\delta-n-1}\
|1-\langle w,a\rangle|^{ r}}}\log^{k}\frac{e}{1-|w|^{2}} $ \ when \ $t-\delta>n+1>r-\delta$.

\vskip2mm

(7) \ \ $\displaystyle{I_{w,a} \asymp\frac{(1-|w|^{2})^{n+1+\delta
-t}}{|1-\langle w,a\rangle|^{r}}\log^{k}\frac{e}{1-|w|^{2}} + \frac{(1-|a|^{2})^{n+1+\delta -r}}{
|1-\langle w,a\rangle|^{ t}}}\log^{k}\frac{e}{1-|a|^{2}}$ \vskip2mm when \ $t-\delta>n+1$, $r-\delta>n+1$.
\vskip2mm
 (8) \ \ $\displaystyle{I_{w,a} \asymp \frac{1}{(1-|w|^{2})^{t-\delta-n-1}|1-\langle w,a\rangle|^{ \delta+n+1}}\log^{k}\frac{e}{1-|w|^{2}}}$
 \vskip2mm + \ $\displaystyle{\frac{1}{|1-\langle
w,a\rangle|^{ t}}\log^{k}\frac{e}{1-|a|^{2}}\log\frac{e}{|1-\langle a,\varphi_{a}(w)\rangle|}}$
 \ when \ $t-\delta>n+1=r-\delta$.
\vskip2mm
 (9) \ \ $\displaystyle{I_{w,a} \asymp \log^{k+1}\frac{e}{1-|w|^{2}}}$ \ when \ $t-\delta=n+1$ \ and \ $r=0$.

 \vskip2mm
{\bf Proposition C} \ Let \ $r+n> 0$ \ and \ $t+n> 0$.  Then  the integral
\begin{eqnarray*}
 &\;& I_{w,a}=\int_{\partial B}\frac{d\sigma(\xi)}{|1-\langle \xi, w\rangle|^{n+ t}\ |1-\langle \xi,a\rangle|^{n+ r}} \ \ (w,a\in B_{n})
 \end{eqnarray*}
has the following bidirectional estimates:
\vskip2mm (1) \ \ $ I_{w,a}\asymp 1 $ \ when \ $t+r+n<0$.
\vskip2mm (2) \ \ $ \displaystyle{I_{w,a}\asymp \log\frac{e}{|1-\langle
w,a\rangle|} }$ \ when \ $t+r+n=0$.
\vskip2mm (3) \ \ $ \displaystyle{I_{w,a}\asymp\frac{1}{|1-\langle
w,a\rangle|^{n+r}}\log\frac{e}{|1-\langle w,\varphi_{w}(a)\rangle|} }
$ \ when \ $t=0>r$.
\vskip2mm (4) \ \ $ \displaystyle{I_{w,a}\asymp\frac{1}{|1-\langle
w,a\rangle|^{t+r+n}} }
$ \ when \ $t+r+n>0>\max\{r,  t\}$.
\vskip2mm
(5) \ \ $ \displaystyle{I_{w,a}\asymp \frac{1}{|1-\langle
w,a\rangle|^{n}}\log\frac{e}{1-|\varphi_{w}(a)|^{2}}} $ \ when \ $t=0=r$.
\vskip2mm
(6) \ \
$ \displaystyle{I_{w,a}\asymp \frac{1}{(1-|w|^{2})^{ t}\
|1-\langle w,a\rangle|^{ n+r}}} $ \ when \ $t>0>r$.

\vskip2mm

(7) \ \ $\displaystyle{I_{w,a} \asymp\frac{(1-|w|^{2})^{
-t}}{|1-\langle w,a\rangle|^{n+ r}} + \frac{(1-|a|^{2})^{ -r}}{
|1-\langle w,a\rangle|^{ n+t}}}$ \ when \ $t>0$, $r>0$.
\vskip2mm
 (8) \ \ $\displaystyle{I_{w,a} \asymp \frac{(1-|w|^{2})^{
-t}}{|1-\langle w,a\rangle|^{ n}} + \frac{1}{|1-\langle
w,a\rangle|^{ n+t}}\log\frac{e}{1-|\varphi_{a}(w)|^{2}}}$
 \ when \ $t>0=r$.
\vskip3mm
{\bf Note 1} \ When $t-\delta>n+1=r-\delta$, it follows from Lemma 2.2 in [10] that
$$
I_{w,a} \asymp \frac{(1-|w|^{2})^{\delta+n+1
-t}}{|1-\langle w,a\rangle|^{ \delta+n+1}} + \frac{1}{|1-\langle
w,a\rangle|^{ t}}\log\frac{e}{1-|\varphi_{a}(w)|^{2}}=L_{1}.
$$
However, it follows from Theorem 3.1 (case $k=0$) in [12]  that
$$
J_{w,a}
\asymp \frac{(1-|w|^{2})^{\delta+n+1
-t}}{|1-\langle w,a\rangle|^{ \delta+n+1}} + \frac{1}{|1-\langle
w,a\rangle|^{ t}}\log\frac{e}{|1-\langle\varphi_{a}(w),a\rangle|}=L_{2}.
$$

Although the expressions of  $L_{1}$ and $L_{2}$ are different,  $L_{1}$ and $L_{2}$  are equivalent.

\vskip2mm
In fact, $L_{2}\lesssim L_{1}$ is obvious. Conversely, let $\displaystyle{\sup_{0<x<2}x^{t-\delta-n-1}\log\frac{e}{x}=M. }$ Then
$$
\frac{ \frac{1}{|1-\langle
w,a\rangle|^{ t}}\log\frac{e|1-\langle w,a\rangle|}{1-|w|^{2}}}{\frac{(1-|w|^{2})^{\delta+n+1
-t}}{|1-\langle w,a\rangle|^{ \delta+n+1}}}=\left(\frac{1-|w|^{2}}{|1-\langle w,a\rangle|}\right)^{t-\delta-n-1}\log\frac{e|1-\langle w,a\rangle|}{1-|w|^{2}}\leq M \ \Rightarrow
$$
 $$
L_{1}\leq L_{2}+\frac{1}{|1-\langle
w,a\rangle|^{ t}}\log\frac{e|1-\langle w,a\rangle|}{1-|w|^{2}}
\leq  L_{2}+\frac{M(1-|w|^{2})^{\delta+n+1
-t}}{|1-\langle w,a\rangle|^{ \delta+n+1}}\leq (M+1)L_{2}.
$$
\vskip3mm
In practical applications, the following integrals are often encountered (for example, Zhou and Chen took advantage of the case $k=2$ in [8]):
\begin{eqnarray*}
 &\;&
G(w)=\displaystyle{\int_{S_{n}}\frac{1}{|1-\langle \xi,w\rangle|^{n+c}}\left|\log\frac{e}{1-\langle \xi,w\rangle}\right|^{k}}\ d\sigma(\xi) \ \ \mbox{and}\\&\;&
F(w)=\displaystyle{\int_{B_{n}}\frac{(1-|z|^{2})^{\delta}}{|1-\langle z,w\rangle|^{n+1+\delta+c}}\left|\log\frac{e}{1-\langle z,w\rangle}\right|^{k}}\ dv(z) \ \ (w\in B_{n}),
\end{eqnarray*}
where $\delta>-1$, \ $c$ and $k$ are real numbers.
\vskip2mm
There is a natural problem. Do  $G(w)$ and $F(w)$ have  bidirectional estimates similar to Proposition A? In this paper, we first discuss this problem.
Since $k$ is an abstract real number, the original proof method of Proposition A is no longer suitable for  $F(w)$ and $G(w)$. Therefore, we need to deal with the two integrals in a completely different way. Otherwise, for two  points in $B_{n}$, we also consider the bidirectional estimates of the following integral in some cases:
$$ J_{w,\eta}=\int_{B_{n}}\frac{(1-|z|^{2})^{
\delta}}{|1-\langle z,w\rangle|^{ t}\ |1-\langle z,\eta\rangle|^{ r}}\log^{-k}\frac{e}{|1-\langle z,\eta\rangle|}\log^{k}\frac{e}{1-|z|^{2}}\
dv(z) \ \ (w, \eta\in B_{n}),$$
where \ $\delta>-1$, \ $r> 0$, \ $t> 0$, \ $k>0$.

\vskip2mm

In the study of function space theory, we often encounter the following problems:
\vskip2mm
(1) \  Let $g$ be a given function.  Suppose that $X$ and $Y$ are two function spaces. Does $fg$  belong to $Y$ if
$f$ belongs to $X$?

\vskip2mm
(2) \ What conditions do the function $g$ need if we want to have $fg\in Y$ for all $f\in X$?
In particular, does $fg$ still belong to $X$ if $f$ and $g$ belong to $X$ ?

\vskip2mm
This is the  pointwise multiplier problem between function spaces. In general, a function in space $X$ multiplied by another function may not be in space $Y$. For example, we take $\psi(w_{1},w_{2})=w_{2}$ and $X=Y=\mathcal{B}(U^{2})$, the Bloch space with 2-dimensional unit cylinder as support domain. It is easy to prove that
 $$
 f(w_{1},w_{2})=\log\frac{1}{1-w_{1}}+\log\frac{1}{1-w_{2}}\in \mathcal{B}(U^{2}).
 $$
 But $\psi f$ does not belong to $\mathcal{B}(U^{2})$. This shows that it is very meaningful for us to consider  pointwise  multipliers  between function spaces.
\vskip2mm
 {\bf  Definition 1.5} \ Let $X$ and $Y$ be two function spaces on domain $\Omega$. If  $\psi f$ belongs to $Y$ for all $f\in X$, then
$\psi$ is called a pointwise  multiplier from  $X$ to $Y$.

\vskip1mm

The study of multiplier theory in function spaces has a long history, and there have been a lot of research results, which are closely related to this paper, such as [5]-[6], [20]-[30] etc. As some applications of these integral estimates, we discuss pointwise multipliers from $F(p,\mu,s)$ to $\mathcal{B_{\nu}}(B_{n})$ in this paper.

\section{\bf  Some Lemmas  \label{sect.2}}

\ \ \ In order to prove the main results, we first give several lemmas.

\vskip2mm { \bf Lemma 2.1 } \ Let \ $c\geq 0$, \ $\delta>-1$, \ $k$ be real. Then the integrals
 \begin{eqnarray*}
 &\;&
 I_{1}(\rho)=\int_{0}^{1}\frac{(1-r)^{\delta}}{(1-r\rho)^{\delta+1+c}}\log^{k}\frac{e}{1-\rho r}\ dr \ \ \mbox{and}\\&\;&
 I_{2}(\rho)=\int_{0}^{1}\frac{(1-r)^{\delta}}{(1-r\rho)^{\delta+1+c}}\log^{k}\frac{e(1-\rho r)}{1-\rho }\ dr \ \ (0\leq \rho<1)
 \end{eqnarray*}
have the following
 bidirectional estimates:
\vskip2mm
(1) \ $I_{1}(\rho)\asymp I_{2}(\rho)\asymp 1$ \ when  $c=0$ and $k<-1$.

\vskip2mm
(2) \  $\displaystyle{ I_{1}(\rho)\asymp  \frac{1}{(1-\rho)^{c}}\log^{k}\frac{e}{1-\rho}}$ \ when $c>0$.
\vskip2mm
(3) \  $\displaystyle{I_{1}(\rho)\asymp I_{2}(\rho)\asymp  \log^{k+1}\frac{e}{1-\rho}}$ \ when  $c=0$ and $k>-1$.
\vskip2mm
(4) \ $\displaystyle{ I_{1}(\rho)\asymp I_{2}(\rho)\asymp  \log\log\frac{e^{2}}{1-\rho}}$ \ when $c=0$ and $k=-1$.
\vskip2mm
{\bf Proof} \ If there exists a constant $0<\rho_{0}<1$ such that $0\leq\rho\leq \rho_{0}$, then it is obvious that these quantities are equivalent.  We only need to consider $\rho$ that is sufficiently close to 1. Therefore, we may let $\rho>1/2$.
\vskip2mm
 By a change of variables  \ $\displaystyle{x=\frac{(1-r)\rho}{1-r\rho}}$, \ we have that
$$
 I_{1}(\rho)=\frac{1}{(1-\rho)^{c}\rho^{\delta+1}}\int_{0}^{\rho}x^{\delta}(1-x)^{c-1}\log^{k}\frac{e(1-x)}{1-\rho}\ dx.
$$

It is clear that $$
\displaystyle{\int_{0}^{\frac{1}{2}}x^{\delta}(1-x)^{c-1}\log^{k}\frac{e(1-x)}{1-\rho}\ dx\asymp \log^{k}\frac{e}{1-\rho}}.$$

 Otherwise, let \ $y=\displaystyle{\frac{1-x}{1-\rho}}$. \ We have that
\begin{eqnarray*}
 &\;&
\int_{\frac{1}{2}}^{\rho}x^{\delta}(1-x)^{c-1}\log^{k}\frac{e(1-x)}{1-\rho}\ dx\asymp \int_{\frac{1}{2}}^{\rho}(1-x)^{c-1}\log^{k}\frac{e(1-x)}{1-\rho}\ dx\\
&\;&
=(1-\rho)^{c}\int_{1}^{\frac{1}{2(1-\rho)}}y^{c-1}\log^{k}ey\ dy.
\end{eqnarray*}

This means that
$$
 I_{1}(\rho)\asymp \frac{1}{(1-\rho)^{c}}\log^{k}\frac{e}{1-\rho}+\int_{1}^{\frac{1}{2(1-\rho)}}y^{c-1}\log^{k}ey\ dy.
$$

By  changes of variables \ $\displaystyle{x=\frac{(1-r)\rho}{1-\rho}}$ \ and \ $y=1+x$, \ we have that
$$
I_{2}(\rho)=\frac{1}{(1-\rho)^{c}\rho^{\delta+1}}\int_{0}^{\frac{\rho}{1-\rho}}\frac{x^{\delta}}{(1+x)^{\delta+1+c}}\log^{k}e(1+x)\ dx
$$
$$
\asymp \frac{1}{(1-\rho)^{c}}\left\{1+\int_{2}^{\frac{1}{1-\rho}}\frac{1}{y^{1+c}}\log^{k}ey\ dy\right\}. \ \ \ \
$$

If \ $\rho\rightarrow 1^{-}$, \ then there are the following results:
\begin{eqnarray*}
 &\;&
\int_{1}^{\frac{1}{2(1-\rho)}}y^{c-1}\log^{k}ey\ dy\asymp \frac{1}{(1-\rho)^{c}}\log^{k}\frac{e}{1-\rho} \ \ \mbox{when \ $c>0$},\\&\;&
\int_{1}^{\frac{1}{2(1-\rho)}}y^{-1}\log^{k}ey\ dy\asymp\int_{2}^{\frac{1}{1-\rho}}\frac{1}{y}\log^{k}ey\ dy\asymp \log^{k+1}\frac{e}{1-\rho} \ \ \mbox{when \ $k>-1$},
\\ &\;&
\int_{1}^{\frac{1}{2(1-\rho)}}y^{-1}\log^{-1}ey\ dy\asymp\int_{2}^{\frac{1}{1-\rho}}\frac{1}{y}\log^{-1}ey\ dy\asymp \log\log\frac{e^{2}}{1-\rho},
\\
&\;&
\int_{1}^{\frac{1}{2(1-\rho)}}y^{-1}\log^{k}ey\ dy\asymp\int_{2}^{\frac{1}{1-\rho}}\frac{1}{y}\log^{k}ey\ dy\asymp 1 \ \ \mbox{when \ $k<-1$}.
\end{eqnarray*}

 This proof is complete. \ \ $\Box$

\vskip2mm
{\bf Lemma 2.2} \ For
 \ $t>n>r> 0$ \ and \ $k<0$, \ let
$$ L_{w,\eta}=\int_{S_{n}}\frac{1}{|1-\langle \xi,w\rangle|^{ t}\ |1-\langle \xi,\eta\rangle|^{ r}}\log^{k}\frac{e}{|1-\langle \xi,\eta\rangle|}\
d\sigma(\xi) \ \ (w, \eta\in B_{n}). \ \ \mbox{Then}$$
$$
L_{w,\eta}\gtrsim \frac{1}{(1-|w|^{2})^{t-n}|1-\langle w,\eta\rangle|^{r}}\log^{k}\frac{e}{|1-\langle w,\eta\rangle|}.
$$

\vskip2mm
{\bf Proof} \ For $\varepsilon>0$ and real number $y$, it is easy to obtain that
$$
1\leq\sup_{0<x<2}\ x^{\varepsilon}\log^{y}\frac{e}{x}\leq \max\left\{e^{\varepsilon-|y|}\left(\frac{|y|+1}{\varepsilon}\right)^{|y|}, \ 2^{\varepsilon}\log^{y}\frac{e}{2}\right\}.\eqno{(2.1)}
$$

For any $u\in S_{n}$ \ and \ $k< 0$,  it follows from (2.1) that
\begin{eqnarray*}
&\;&
|1-\langle u,\varphi_{ w}(\eta)\rangle|^{r}\log^{-k}\frac{e|1-\langle u,w\rangle|}{|1-\langle w,\eta\rangle||1-\langle u,\varphi_{ w}(\eta)\rangle|}\\
&\;&
\lesssim|1-\langle u,\varphi_{ w}(\eta)\rangle|^{r}\left\{\log^{-k}\frac{e}{|1-\langle w,\eta\rangle|}+\log^{-k}\frac{e}{|1-\langle u,\varphi_{ w}(\eta)\rangle|}\right\}\\
&\;&
\lesssim \log^{-k}\frac{e}{|1-\langle w,\eta\rangle|}+1\asymp\log^{-k}\frac{e}{|1-\langle w,\eta\rangle|}.
\end{eqnarray*}

 By a change of variables  $\xi=\varphi_{w}(u)$, (1.1)-(1.2), (4.7) in [1], the increasing property of integral mean of holomorphic function, we obtain that
\begin{eqnarray*}
&\;&
L_{w,\eta}=\frac{(1-|w|^{2})^{n-t}}{|1-\langle w,\eta\rangle|^{r}}\int_{S_{n}}\frac{|1-\langle u, w\rangle|^{t+r-2n} }{|1-\langle u,\varphi_{w}(\eta)\rangle|^{r}}\log^{k}\frac{e|1-\langle u,w\rangle|}{|1-\langle w,\eta\rangle||1-\langle u,\varphi_{ w}(\eta)\rangle|}\ d\sigma(u)
\\&\;&
\gtrsim\frac{1}{(1-|w|^{2})^{t-n}|1-\langle w,\eta\rangle|^{r}}\log^{k}\frac{e}{|1-\langle w,\eta\rangle|}\int_{S_{n}}|1-\langle u, w\rangle|^{t+r-2n}\ d\sigma(u)
\\&\;&
\geq\frac{1}{(1-|w|^{2})^{t-n}|1-\langle w,\eta\rangle|^{r}}\log^{k}\frac{e}{|1-\langle w,\eta\rangle|}.
\end{eqnarray*}

This proof is complete. \ \ $\Box$

\vskip2mm
 {  \bf Lemma 2.3 ([31])} \ Let $\mu$ be a normal function on $[0,1)$. Suppose that $a$ and $b$ are the two parameters in the definition of $\mu$. Then there are the following results.
 \vskip2mm
 (1) \ \ For all \ $z,w\in B_{n}$,
 $$
 \frac{\mu(|z|)}{\mu(|w|)}\leq \left(\frac{1-|z|^{2}}{1-|w|^{2}}\right)^{a}+\left(\frac{1-|z|^{2}}{1-|w|^{2}}\right)^{b}.
 $$

 (2) \ \ If \ $z\in B_{n}$, \ then  $\mu(|w|)\asymp \mu(|z|)$ for all $w\in D(z,1)$.
\vskip2mm
This Lemma comes from Lemma 2.2 in [31].

 \vskip2mm { \bf Lemma 2.4 } \ Let $\mu$ be a normal function on $[0,1)$, and let $\mu_{1}(r)=(1-r^{2})^{\frac{n-s}{p}}\mu(r)$ be also a normal function on $[0,1)$.  If \ $f\in F(p,\mu,s)$, then
 \begin{eqnarray*}
 &\;&
 |\nabla f(w)|\lesssim \frac{||f||_{p,\mu,s}}{(1-|w|^{2})^{\frac{n-s}{p}}\mu(|w|)} \ \ \mbox{and}\\
 &\;&
 |f(w)|\lesssim \left\{1+\int_{0}^{|w|}\frac{d\rho}{(1-\rho^{2})^{\frac{n-s}{p}}\mu(\rho)}\right\}||f||_{p,\mu,s} \ \ \mbox{for all $w\in B_{n}$.}
 \end{eqnarray*}
 In particular, $F(p,\mu,s)=\mathcal{B}_{\mu_{1}}(B_{n})$ when $s>n$.

\vskip2mm
{\bf Proof} \ For any \ $l\in \{1,2,\cdots,n\}$ and  $w\in B_{n}$, by  Lemma 2.20 and Lemma 2.24 in [1], Lemma 2.3, we have that
\begin{eqnarray*}
 &\;&
 \left|\frac{\partial f}{\partial w_{l}}(w)\right|^{p}\lesssim \frac{1}{(1-|w|^{2})^{n+1}}\int_{D(w,1)}\left|\frac{\partial f}{\partial z_{l}}(z)\right|^{p}\ dv(z)\\
 &\;&
 \lesssim \frac{(1-|w|^{2})^{s-n}}{\mu^{p}(|w|)}\int_{D(w,1)}\frac{(1-|w|^{2})^{s}|\nabla f(z)|^{p}\ \mu^{p}(|z|)}{|1-\langle z,w\rangle|^{2s}(1-|z|^{2})}\ dv(z)\\
 &\;&
 \leq \frac{(1-|w|^{2})^{s-n}}{\mu^{p}(|w|)}\ ||f||_{p,\mu,s}^{p}.
 \end{eqnarray*}
  This shows that
  $$|\nabla f(w)|\asymp \sum_{l=1}^{n} \left|\frac{\partial f}{\partial w_{l}}(w)\right|\asymp \left\{\sum_{l=1}^{n} \left|\frac{\partial f}{\partial w_{l}}(w)\right|^{p}\right\}^{\frac{1}{p}}\lesssim \frac{||f||_{p,\mu,s}}{(1-|w|^{2})^{\frac{n-s}{p}}\mu(|w|)}.
$$
Therefore, $ F(p,\mu,s)\subseteq \mathcal{B}_{\mu_{1}}(B_{n})$.
Moreover, we have that
$$
|f(w)|=\left|f(0,\cdots,0)+\int_{0}^{1}\langle \nabla f(\rho w), \overline{w}\rangle\ d\rho\right|$$$$\lesssim  \left\{1+\int_{0}^{|w|}\frac{d\rho}{(1-\rho^{2})^{\frac{n-s}{p}}\mu(\rho)}\right\}||f||_{p,\mu,s}.
$$

 When $s>n$, for any $f\in \mathcal{B}_{\mu_{1}}(B_{n})$ and $w\in
B_{n}$, it follows from Proposition A that
\begin{eqnarray*}
 &\;&
 \int_{B_{n}}\frac{(1-|w|^{2})^{s}|\nabla f(z)|^{p}\ \mu^{p}(|z|)}{|1-\langle z,w\rangle|^{2s}(1-|z|^{2})}\ dv(z)\\
 &\;&\lesssim ||f||_{ \mathcal{B}_{\mu_{1}}}^{p}\int_{B_{n}}\frac{(1-|w|^{2})^{s}(1-|z|^{2})^{s-n-1}}{|1-\langle z,w\rangle|^{2s}}\ dv(z)\\
 &\;&
 \asymp ||f||_{ \mathcal{B}_{\mu_{1}}}^{p} \ \Rightarrow f\in F(p,\mu,s) \ \Rightarrow \mathcal{B}_{\mu_{1}}(B_{n})\subseteq F(p,\mu,s).
\end{eqnarray*}

This proof is complete. \ \ $\Box$
\vskip2mm
{ \bf Lemma 2.5 ([32, 33])}  \ Let $\mu$ be a normal function on $[0,1)$ and
$$
g(u)=1+\sum_{j=1}^{\infty}2^{j}\ u^{n_{j}}\ \ (u\in D),
$$
where $n_{j}$ is the integral part of $(1-r_{j})^{-1}$, \ $\mu(r_{j})=2^{-j}$ \ ($j=1,2,...$).
Then there are the following results:
\vskip2mm
(1) \ $g(r)$ is increasing on $[0,1)$, and there exist $M_{0}>0$ and $M_{1}<\infty$ such that
$$
\displaystyle{\inf_{r\in [0,1)}\mu(r)g(r)\geq M_{0}}, \ \ \displaystyle{\sup_{u\in
D}\mu(|u|)|g(u)|\leq M_{1}}$$
\vskip2mm
(2) \ There exists a constant $M_{2}<\infty$ such that
$$\displaystyle{\sup_{r\in [0,1)}(1-r)\mu(r)g'(r)\leq M_{2}}.$$
\vskip2mm
The first part of this Lemma comes from Theorem 1 in [32], and the second part comes from Lemma 2.4 in [33].
\vskip 2mm
 {\bf Lemma 2.6 ([1])} \ \ For \ $\alpha>-1$, \ if \ $f\in A_{\alpha}^{1}(B_{n})$, \ then
 $$
 f(z)=\int_{B}\frac{f(w)}{(1-\langle z,w\rangle)^{n+1+\alpha}}\ dv_{\alpha}(w)\ \ (z\in B_{n}).
 $$

This Lemma comes from Theorem 2.2 in [1].

\section{\bf  Main Results  \label{sect.3}}

\ \ \ We first prove  three propositions.
\vskip2mm
{\bf Proposition 3.1} \ Let $c$ and $k$ be real, $\delta>-1$. \ Then the integrals
\begin{eqnarray*}
 &\;&
G(w)=\displaystyle{\int_{S_{n}}\frac{1}{|1-\langle \xi,w\rangle|^{n+c}}\left|\log\frac{e}{1-\langle \xi,w\rangle}\right|^{k}}\ d\sigma(\xi) \ \ \mbox{and}\\&\;&
F(w)=\displaystyle{\int_{B_{n}}\frac{(1-|z|^{2})^{\delta}}{|1-\langle z,w\rangle|^{n+1+\delta+c}}\left|\log\frac{e}{1-\langle z,w\rangle}\right|^{k}}\ dv(z) \ \ (w\in B_{n})
\end{eqnarray*}
have the following bidirectional estimates:

\vskip2mm
(1) \ $G(w)\asymp F(w)\asymp 1$ \ when $c<0$, or $c=0$ and $k<-1$.
\vskip2mm
(2) \ $\displaystyle{G(w)\asymp F(w)\asymp \frac{1}{(1-|w|^{2})^{c}}\log^{k}\frac{e}{1-|w|^{2}} }$ \ when $c>0$.
\vskip2mm
(3) \ $\displaystyle{G(w)\asymp F(w)\asymp \log^{k+1}\frac{e}{1-|w|^{2}} }$ \ when $c=0$ and $k>-1$.
\vskip2mm
(4) \ $\displaystyle{G(w)\asymp F(w)\asymp \log\log\frac{e^{2}}{1-|w|^{2}} }$ \ when $c=0$ and $k=-1$.

\vskip2mm
{\bf Proof} \  If there exists a constant $0<\rho_{0}<1$ such that $1-|w|^{2}\geq \rho_{0}$, then these bidirectional estimates are obvious.  Therefore, we  let $1-|w|^{2}$  be sufficiently close to 0.
When \ $|z|\geq 1$, \ it is clear that $$ \log|z|\leq |\log z|=\sqrt{(\log |z|)^{2}+(\arg z)^{2}}\leq \log e^{\pi}|z|.\eqno{(3.1)}$$
 This shows that
$$
G(w)\asymp\int_{S_{n}}\frac{1}{|1-\langle \xi,w\rangle|^{n+c}}\log^{k}\frac{e}{|1-\langle \xi,w\rangle|}\ d\sigma(\xi).
$$

When  $c<0$, it follows from (2.1) that we may take \ $c< c'<0$ \ such that
$$
 \frac{1}{|1-\langle \xi,w\rangle|^{n+c}}\log^{k}\frac{e}{|1-\langle \xi,w\rangle|}\lesssim \frac{1}{|1-\langle \xi,w\rangle|^{n+c'}}
$$
for all $\xi\in S_{n}$ and $w\in B_{n}$. It follows from the increasing property of integral mean of holomorphic function and Proposition A  that
$$
1\leq G(w)\lesssim \int_{S_{n}}\frac{d\sigma(\xi)}{|1-\langle \xi,w\rangle|^{n+c'}}\asymp 1.
$$

 By  a change of variables $\xi=\varphi_{w}(\eta)$, (4.7) in [1],  (1.1)-(1.2), we have that
\begin{eqnarray*}
 &\;&
 G(w)\asymp\frac{1}{(1-|w|^{2})^{c}}\int_{S_{n}}\frac{1}{|1-\langle \eta,w\rangle|^{n-c}}\log^{k}\frac{e|1-\langle \eta,w\rangle|}{1-|w|^{2}}\ d\sigma(\eta)=\frac{J(w)}{(1-|w|^{2})^{c}}.
\end{eqnarray*}

Next, we consider $J(w)$ when $c\geq 0$.
\vskip2mm When $n=1$,  it follows from the rotation invariance of the integral   that
\begin{eqnarray*}
 &\;&
 J(w)=\int_{-\pi}^{\pi}\frac{1}{|1-|w|e^{i\theta}|^{1-c}}\log^{k}\frac{e|1-|w|e^{i\theta}|}{1-|w|^{2}}\ \frac{d\theta}{2\pi}\\
 &\;&
 =\int_{0}^{\pi}\frac{1}{(1+|w|^{2}-2|w|\cos\theta)^{\frac{1-c}{2}}}\log^{k}\frac{e^{2}(1+|w|^{2}-2|w|\cos\theta)}{(1-|w|^{2})^{2}}\ \frac{d\theta}{2^{k}\pi}\\
 &\;&
 =\frac{1}{2^{k}\pi} \int_{-1}^{1}\frac{(1-x^{2})^{-\frac{1}{2}}}{(1+|w|^{2}-2|w|x)^{\frac{1-c}{2}}}\log^{k}\frac{e^{2}(1+|w|^{2}-2|w|x)}{(1-|w|^{2})^{2}}\ dx\\
 &\;&
 \asymp \log^{k}\frac{e}{1-|w|^{2}}+\int_{0}^{1}\frac{(1-x)^{-\frac{1}{2}}}{(1+|w|^{2}-2|w|x)^{\frac{1-c}{2}}}\log^{k}\frac{e^{2}(1+|w|^{2}-2|w|x)}{(1-|w|^{2})^{2}}\ dx.
\end{eqnarray*}

Without losing generality, we let $|w|>1/2$. By   a change of variables  $$\rho=\frac{(1+|w|^{2})(1-x)}{1+|w|^{2}-2|w|x}, \ \ \mbox{we have that}$$
\begin{eqnarray*}
 J(w)&\asymp& \log^{k}\frac{e}{1-|w|^{2}}+(1-|w|)^{c}\int_{0}^{1}\frac{\log^{k}\frac{e}{1-\frac{2|w|\rho}{1+|w|^{2}}}}{\rho^{\frac{1}{2}}\left(1-\frac{2|w|\rho}{1+|w|^{2}}\right)^{1+\frac{c}{2}}}\ d\rho\\
 &\asymp& \log^{k}\frac{e}{1-|w|^{2}}+(1-|w|)^{c}\left\{1+\int_{\frac{1}{2}}^{1}\frac{\log^{k}\frac{e}{1-\frac{2|w|\rho}{1+|w|^{2}}}}{\left(1-\frac{2|w|\rho}{1+|w|^{2}}\right)^{1+\frac{c}{2}}}\ d\rho\right\}
\end{eqnarray*}
$$
\asymp  \ \ \log^{k}\frac{e}{1-|w|^{2}}+(1-|w|)^{c}\int_{0}^{1}\frac{\log^{k}\frac{e}{1-\frac{2|w|\rho}{1+|w|^{2}}}}{\left(1-\frac{2|w|\rho}{1+|w|^{2}}\right)^{1+\frac{c}{2}}}\ d\rho.
$$
It follows from Lemma 2.1 that we can get the estimates of $J(w)$ by different cases.
\vskip2mm
When $n>1$, it follows from (1.13) in [1] that
\begin{eqnarray*}
 &\;&
 J(w)=(n-1)\int_{D}\frac{(1-|z|^{2})^{n-2}}{|1-|w|z|^{n-c}}\log^{k}\frac{e|1-|w|z|}{1-|w|^{2}}\ dA(z).
\end{eqnarray*}

If $c>0$ and $k\geq 0$, then it follows from Proposition A that
$$
J(w)\lesssim \log^{k}\frac{e}{1-|w|^{2}}\int_{D}\frac{(1-|z|^{2})^{n-2}}{|1-|w|z|^{n-c}}\ dA(z)\asymp \log^{k}\frac{e}{1-|w|^{2}}.
$$

If $c>0$ and $k< 0$, then we take $0<\varepsilon<c$.  It is easy to obtain that
$$
\sup_{1-|w|\leq x\leq 1+|w|}x^{\varepsilon}\log^{k}\frac{ex}{1-|w|^{2}}\asymp \max\left\{(1-|w|^{2})^{\varepsilon}, \log^{k}\frac{e}{1-|w|}\right\}\lesssim \log^{k}\frac{e}{1-|w|^{2}}.
$$

Therefore, it follows from Proposition A that
$$
J(w)\lesssim \log^{k}\frac{e}{1-|w|^{2}}\int_{D}\frac{(1-|z|^{2})^{n-2}}{|1-|w|z|^{n-c+\varepsilon}}\ dA(z)\asymp \log^{k}\frac{e}{1-|w|^{2}}.
$$

On the other hand, we have that
$$
J(w) \gtrsim \int_{|z|\leq\frac{1}{2}}\frac{(1-|z|^{2})^{n-2}}{|1-|w|z|^{n-c}}\log^{k}\frac{e|1-|w|z|}{1-|w|^{2}}\ dA(z)\asymp \log^{k}\frac{e}{1-|w|^{2}}.
$$

This means that  $\displaystyle{J(w)\asymp\log^{k}\frac{e}{1-|w|^{2}}}$ \ when $c>0$. Therefore,
$$\displaystyle{G(w)\asymp\frac{1}{(1-|w|^{2})^{c}}\log^{k}\frac{e}{1-|w|^{2}}} \ \ \mbox{when $c>0$}.$$

For any $1/2<\rho<1$ and any real number $k$,  let $\displaystyle{x=\frac{2\rho|w|(1-r)}{(1-\rho|w|)^{2}}}$.  Similar to the previous calculation, we can obtain that
\begin{eqnarray*}
 &\;&
\int_{-\pi}^{\pi}\frac{1}{|1-\rho|w|e^{i\theta}|^{n}}\log^{k}\frac{e|1-\rho|w|e^{i\theta}|}{1-|w|^{2}}\ \frac{d\theta}{2\pi}\asymp \log^{k}\frac{e}{1-|w|^{2}}
\\
 &\;&
  + \ \int_{0}^{1}\frac{1}{(1+\rho^{2}|w|^{2}-2\rho|w|r)^{\frac{n}{2}}}\log^{k}\frac{e^{2}(1+\rho^{2}|w|^{2}-2\rho|w|r)}{(1-|w|^{2})^{2}}\frac{ dr}{\sqrt{1-r}}
 \\ &\;& \asymp \log^{k}\frac{e}{1-|w|^{2}}
 + \ \frac{1}{(1-\rho|w|)^{n-1}}\int_{0}^{\frac{2\rho|w|}{(1-\rho |w|)^{2}}}\frac{\log^{k}\frac{e^{2}(1-\rho|w|)^{2}}{(1-|w|^{2})^{2}}(1+x)}{x^{\frac{1}{2}}(1+x)^{\frac{n}{2}}}\ dx.
\end{eqnarray*}

It is clear that
$$
\int_{0}^{\frac{8}{9}}\frac{1}{x^{\frac{1}{2}}(1+x)^{\frac{n}{2}}}\log^{k}\frac{e^{2}(1-\rho|w|)^{2}}{(1-|w|^{2})^{2}}(1+x)\ dx\asymp \log^{k}\frac{e(1-\rho|w|)}{1-|w|}.
$$

When  $k\geq 0$, we have that $$\displaystyle{\log^{k}\frac{e^{2}(1-\rho|w|)^{2}}{(1-|w|^{2})^{2}}(1+x)\asymp \log^{k}\frac{e^{2}(1-\rho|w|)^{2}}{(1-|w|^{2})^{2}}+\log^{k}(1+x)}. \ \mbox{ Therefore,}$$
\begin{eqnarray*}
 &\;&
\int_{\frac{8}{9}}^{\frac{2\rho|w|}{(1-\rho |w|)^{2}}}\frac{1}{x^{\frac{1}{2}}(1+x)^{\frac{n}{2}}}\log^{k}\frac{e^{2}(1-\rho|w|)^{2}}{(1-|w|^{2})^{2}}(1+x)\ dx\\
&\;&
 \lesssim\int_{\frac{8}{9}}^{\infty}\frac{1}{x^{\frac{n+1}{2}}}\left\{\log^{k}\frac{e(1-\rho|w|)}{1-|w|^{2}}+\log^{k}(x+1)\right\}\ dx\\
 &\;&
 \asymp \log^{k}\frac{e(1-\rho|w|)}{1-|w|}+1\asymp\log^{k}\frac{e(1-\rho|w|)}{1-|w|}.
\end{eqnarray*}

When $k<0$, we have that $$\displaystyle{\log^{k}\frac{e^{2}(1-\rho|w|)^{2}}{(1-|w|^{2})^{2}}(1+x)\leq \log^{k}\frac{e^{2}(1-\rho|w|)^{2}}{(1-|w|^{2})^{2}}}. \ \mbox{ Therefore,}$$
\begin{eqnarray*}
 &\;&
\int_{\frac{8}{9}}^{\frac{2\rho|w|}{(1-\rho |w|)^{2}}}\frac{1}{x^{\frac{1}{2}}(1+x)^{\frac{n}{2}}}\log^{k}\frac{e^{2}(1-\rho|w|)^{2}}{(1-|w|^{2})^{2}}(1+x)\ dx\\
&\;&
\lesssim \log^{k}\frac{e(1-\rho|w|)}{1-|w|^{2}}\int_{\frac{8}{9}}^{\infty}\frac{1}{x^{\frac{n+1}{2}}}\ dx\asymp\log^{k}\frac{e(1-\rho|w|)}{1-|w|}.
\end{eqnarray*}

This means that
\begin{eqnarray*}
 &\;&
\int_{-\pi}^{\pi}\frac{1}{|1-\rho|w|e^{i\theta}|^{n}}\log^{k}\frac{e|1-\rho|w|e^{i\theta}|}{1-|w|}\ \frac{d\theta}{2\pi}
\\
&\;&\asymp \log^{k}\frac{e}{1-|w|^{2}}+\frac{1}{(1-\rho|w|)^{n-1}}\log^{k}\frac{e(1-\rho|w|)}{1-|w|}.
\end{eqnarray*}

If $c=0$, then it follows from the polar coordinate and the above result that
\begin{eqnarray*}
 &\;&
J(w)\asymp \log^{k}\frac{e}{1-|w|^{2}}\int_{0}^{1}(1-\rho)^{n-2}\ d\rho+\int_{0}^{1}\frac{(1-\rho)^{n-2}}{(1-\rho|w|)^{n-1}}\log^{k}\frac{e(1-\rho|w|)}{1-|w|}\ d\rho.
\end{eqnarray*}

It follows from Lemma 2.1 that  we can get the estimates of $J(w)$ for all cases.

\vskip2mm
Finally, we consider $F(w)$.
\vskip2mm

When $c<0$, let \ $c< c'<0$. By (3.1), (2.1) and Proposition A, we have that
$$
1\lesssim F(w)\asymp\int_{B_{n}}\frac{(1-|z|^{2})^{\delta}}{|1-\langle z,w\rangle|^{n+1+\delta+c}}\log^{k}\frac{e}{|1-\langle z,w\rangle|}\ dv(z)
$$
$$
\lesssim \int_{B_{n}}\frac{(1-|z|^{2})^{\delta}}{|1-\langle z,w\rangle|^{n+1+\delta+c'}}\ dv(z)\lesssim 1.
$$

When $c\geq 0$, by (3.1), Lemma 1.8 in [1] and the estimate of $G(w)$, we have that
$$
F(w)\asymp \int_{0}^{1}\frac{(1-\rho)^{\delta}}{(1-\rho|w|)^{\delta+1+c}}\log^{k}\frac{e}{1-\rho|w|}\ d\rho,
$$

By Lemma 2.1, we can get  the estimates of $F(w)$ in different cases.
\vskip2mm
The proof is complete. \ \ $\Box$

\vskip2mm
{\bf Proposition 3.2} \ For \ $\delta>-1$,
 \ $r> 0$, \ $t> 0$, \ $k>0$, \ the integral
$$ J_{w,\eta}=\int_{B_{n}}\frac{(1-|z|^{2})^{
\delta}}{|1-\langle z,w\rangle|^{ t}\ |1-\langle z,\eta\rangle|^{ r}}\left|\log\frac{e}{1-\langle z,\eta\rangle}\right|^{-k}\log^{k}\frac{e}{1-|z|^{2}}\
dv(z) $$ $(w, \eta\in B_{n})$ has the following bidirectional estimates:
\vskip2mm
(1) \ When \ $r+t-\delta>n+1> \max\{t-\delta, \ r-\delta\}$,
 $$
\displaystyle{J_{w,\eta}\asymp \frac{1}{|1-\langle w,\eta\rangle|^{r+t-\delta-n-1}}}.$$

(2) \ When \ $t-\delta=n+1>r-\delta$,
$$
J_{w,\eta}\asymp \frac{1}{|1-\langle w,\eta\rangle|^{r}}\log^{-k}\frac{e}{|1-\langle w,\eta\rangle|}\log^{k}\frac{e}{1-|w|^{2}}\log\frac{e}{|1-\langle w,\varphi_{w}(\eta)\rangle|}.
$$

(3) \ When \ $t-\delta>n+1>r-\delta$,
$$
J_{w,\eta}\asymp \frac{1}{(1-|w|^{2})^{t-\delta-n-1}|1-\langle w,\eta\rangle|^{r}}\log^{-k}\frac{e}{|1-\langle w,\eta\rangle|}\log^{k}\frac{e}{1-|w|^{2}}.
$$

{\bf Proof }  Without losing generality, let $1-|w|^{2}$ and $|1-\langle w,\eta\rangle|$ be sufficiently close to 0. It follows from (3.1) that
$$
J_{w,\eta}\asymp\int_{B_{n}}\frac{(1-|z|^{2})^{
\delta}}{|1-\langle z,w\rangle|^{ t}\ |1-\langle z,\eta\rangle|^{ r}}\log^{-k}\frac{e}{|1-\langle z,\eta\rangle|}\log^{k}\frac{e}{1-|z|^{2}}\
dv(z).
$$

By a partition of $B_{n}$, there exist
$\Omega_{1}$, $\Omega_{2}$, $\Omega_{3}$,
$\Omega_{4}$ such that $B=\Omega_{1}\cup\Omega_{2}\cup\Omega_{3}\cup\Omega_{4}$, where $\Omega_{j}$ and $\Omega_{k}$ \ ($j\neq k$) are mutually disjoint. By Lemma 3.3 in [34], we have
\begin{eqnarray*}
&\;&
|1-\langle z,w\rangle|\lesssim|1-\langle w,\eta\rangle|\lesssim|1-\langle z,\eta\rangle| \ \ (z\in \Omega_{1});\\
&\;&|1-\langle z,\eta\rangle|\lesssim|1-\langle w,\eta\rangle|\lesssim|1-\langle z,w\rangle| \ \ (z\in \Omega_{2});
\\
&\;&|1-\langle w,\eta\rangle|\lesssim|1-\langle z,w\rangle|\lesssim|1-\langle z,\eta\rangle| \ \ (z\in \Omega_{3});\\
&\;&|1-\langle w,\eta\rangle|\lesssim|1-\langle z,\eta\rangle|\lesssim|1-\langle z,w\rangle| \ \ (z\in \Omega_{4}).
\end{eqnarray*}

When $z\in \Omega_{1}\cup \Omega_{3}$, we have that $|1-\langle z,\eta\rangle|\gtrsim |1-\langle w,\eta\rangle|$. If \ $t-\delta>n+1$, then it follows from Proposition B(6) that
$$
J_{1}=\int_{\Omega_{1}\cup \Omega_{3}}\frac{(1-|z|^{2})^{
\delta}}{|1-\langle z,w\rangle|^{ t}\ |1-\langle z,\eta\rangle|^{ r}}\log^{-k}\frac{e}{|1-\langle z,\eta\rangle|}\log^{k}\frac{e}{1-|z|^{2}}\
dv(z)
$$
$$
\lesssim \frac{1}{|1-\langle w,\eta\rangle|^{ r}}\log^{-k}\frac{e}{|1-\langle w,\eta\rangle|}\int_{B_{n}}\frac{(1-|z|^{2})^{
\delta}}{|1-\langle z,w\rangle|^{ t}}\log^{k}\frac{e}{1-|z|^{2}}\
dv(z)
$$$$
 \asymp\frac{1}{(1-|w|^{2})^{t-\delta-n-1}|1-\langle w,\eta\rangle|^{ r}}\log^{-k}\frac{e}{|1-\langle w,\eta\rangle|}\log^{k}\frac{e}{1-|w|^{2}}. \ \ \ \ \ \ \eqno{(3.2)}
$$

When $t-\delta=n+1$, we may take $0<\varepsilon<r$. This means that $(r-\varepsilon)-\delta<n+1$ and $r-\varepsilon>0$. By  a change of variables $z=\varphi_{w}(u)$, (1.1)-(1.2),  Proposition B(2), we have that
\begin{eqnarray*}
&\;&
J_{1}\lesssim \frac{1}{|1-\langle w,\eta\rangle|^{\varepsilon}}\log^{-k}\frac{e}{|1-\langle w,\eta\rangle|}\int_{B_{n}}\frac{(1-|z|^{2})^{
\delta}\log^{k}\frac{e}{1-|z|^{2}}\
dv(z)}{|1-\langle z,w\rangle|^{\delta+n+1}|1-\langle z,\eta\rangle|^{r-\varepsilon}}\\
&\;&
=\frac{1}{|1-\langle w,\eta\rangle|^{r}}\log^{-k}\frac{e}{|1-\langle w,\eta\rangle|}\int_{B_{n}}\frac{(1-|u|^{2})^{
\delta}\log^{k}\frac{e}{1-|\varphi_{w}(u)|^{2}}\
dv(u)}{|1-\langle u,w\rangle|^{\delta+n+1-(r-\varepsilon)}|1-\langle u,\varphi_{w}(\eta)\rangle|^{r-\varepsilon}}\\
&\;&
\lesssim \frac{1}{|1-\langle w,\eta\rangle|^{r}}\log^{-k}\frac{e}{|1-\langle w,\eta\rangle|}\int_{B_{n}}\frac{(1-|u|^{2})^{
\delta}\left\{\log^{k}\frac{e}{1-|w|^{2}}+\log^{k}\frac{e}{1-|u|^{2}}\right\}
dv(u)}{|1-\langle u,w\rangle|^{\delta+n+1-(r-\varepsilon)}|1-\langle u,\varphi_{w}(\eta)\rangle|^{r-\varepsilon}}\\
&\;&
\asymp \frac{\log^{-k}\frac{e}{|1-\langle w,\eta\rangle|}}{|1-\langle w,\eta\rangle|^{r}}\left\{\log^{k}\frac{e}{1-|w|^{2}}\log\frac{e}{|1-\langle w,\varphi_{w}(\eta)\rangle|}+\log^{k+1}\frac{e}{|1-\langle w,\varphi_{w}(\eta)\rangle|}\right\} \ \
\end{eqnarray*}
$$
\asymp\frac{1}{|1-\langle w,\eta\rangle|^{r}}\log^{-k}\frac{e}{|1-\langle w,\eta\rangle|}\log^{k}\frac{e}{1-|w|^{2}}\log\frac{e}{|1-\langle w,\varphi_{w}(\eta)\rangle|}. \ \ \ \ \ \eqno{(3.3)}
$$

\vskip2mm

Let $m=|1-\langle w,\eta\rangle|$. When $z\in \Omega_{2}\cup \Omega_{4}$, it is clear that $$|1-\langle z,w\rangle|\gtrsim (m+|1-\langle z,\eta\rangle|)\geq |m+1-\langle z,\eta\rangle|.$$

When $n>1$, it is similar to the proof of (3.17) in [12] (we omit some processes).
If \ $r-\delta<n+1<t+r-\delta$, then it follows from Lemma 1.8 and (1.13) in [1] that
\begin{eqnarray*}
&\;&
J_{2}=\int_{\Omega_{2}\cup \Omega_{4}}\frac{(1-|z|^{2})^{
\delta}}{|1-\langle z,w\rangle|^{ t}\ |1-\langle z,\eta\rangle|^{ r}}\log^{-k}\frac{e}{|1-\langle z,\eta\rangle|}\log^{k}\frac{e}{1-|z|^{2}}\
dv(z)\\
&\;&\lesssim\int_{0}^{1}(1-\rho^{2})^{\delta}\log^{k}\frac{e}{1-\rho^{2}}\left\{\int_{D}\frac{\rho^{2n-1}(1-|\xi|^{2})^{n-2}}{|m+1-\rho\xi|^{t}|1-\rho\xi|^{r}}\log^{-k}\frac{e}{|1-\rho\xi|}\ dA(\xi)\right\}d\rho
\\
&\;&\lesssim \int_{0}^{1}(1-\rho)^{\delta}\log^{k}\frac{e}{1-\rho}\left\{\int_{|u-1|<\rho}\frac{(\rho-|u-1|)^{n-2}}{|m+u|^{t}|u|^{r}}\log^{-k}\frac{e}{|u|}\ dA(u)\right\}d\rho \
\end{eqnarray*}
\begin{eqnarray*}
&\;&
\lesssim \int_{0}^{1}(1-\rho)^{\delta}\log^{k}\frac{e}{1-\rho}\left\{\int_{1-\rho}^{1}\int_{-\frac{\pi}{2}}^{\frac{\pi}{2}}\frac{\{R-(1-\rho)\}^{n-2}}{(m+R)^{t}R^{r-1}}\log^{-k}\frac{e}{R}\ dRd\theta\right\}d\rho \ \ \ \ \ \ \ \
\\
&\;&+ \
\int_{0}^{1}(1-\rho)^{\delta}\log^{k}\frac{e}{1-\rho}\left\{\int_{1}^{1+\rho}\int_{-\frac{\pi}{2}}^{\frac{\pi}{2}}\frac{\rho^{n-2}}{(m+R)^{t}R^{r-1}}\log^{-k}\frac{e}{R}\ dRd\theta\right\}d\rho\\
&\;&\lesssim\int_{0}^{1}(1-\rho)^{\delta}\log^{k}\frac{e}{1-\rho}\left\{\int_{1-\rho}^{1}\frac{R^{n-1-r}}{(m+R)^{t}}\log^{-k}\frac{e}{R}\ dR\right\}d\rho+1
\\
&\;&\lesssim \int_{0}^{1}\frac{(1-\rho)^{\delta}}{m^{r+t-n}}\log^{k}\frac{e}{1-\rho}\left\{\int_{\frac{1-\rho}{m}}^{\frac{1}{m}}\frac{x^{n-1-r}}{(1+x)^{t}}\log^{-k}\frac{e}{mx}\ dx\right\}d\rho +1
\\
&\;&=
\frac{1}{m^{r+t-n}}\int_{0}^{\frac{1}{m}}\frac{x^{n-1-r}}{(1+x)^{t}}\log^{-k}\frac{e}{mx}\left\{\int_{1-mx}^{1}(1-\rho)^{\delta}\log^{k}\frac{e}{1-\rho}\ d\rho\right\}dx +1 \ \ \ \ \ \
 \\
&\;&\asymp \frac{1}{m^{r+t-\delta-n-1}}\int_{0}^{\frac{1}{m}}\frac{x^{n+\delta-r}}{(1+x)^{t}}\ dx+1
\end{eqnarray*}
$$
\leq \frac{1}{m^{r+t-\delta-n-1}}\int_{0}^{1}y^{n+\delta-r}(1-y)^{r+t-\delta-n-2}\ dy\asymp \frac{1}{|1-\langle w,\eta\rangle|^{r+t-\delta-n-1}}.  \eqno{(3.4)}
$$

When $n>1$, it is similar to the proof of (3.4). If \ $t-\delta<n+1<r+t-\delta$, then we can also get that
\begin{eqnarray*}
 J_{1}
 &\lesssim& \int_{\Omega_{1}\cup\Omega_{3}}\frac{(1-|z|^{2})^{
\delta}\log^{k}\frac{e}{1-|z|^{2}}}{|1-\langle z,w\rangle|^{ t}\ |m+1-\langle z,w\rangle|^{ r}}\log^{-k}\frac{e}{|m+1-\langle z,w\rangle|}\
dv(z)
\\ &\lesssim& \frac{1}{m^{r+t-\delta-n-1}}\int_{0}^{\frac{1}{m}}\frac{x^{n+\delta-t}}{(1+x)^{r}\log^{k}\frac{e}{m(1+x)}}\log^{k}\frac{e}{mx}\ dx+1.
\end{eqnarray*}

It is clear that
$$
\int_{0}^{\frac{1}{2}}\frac{x^{n+\delta-t}}{(1+x)^{r}\log^{k}\frac{e}{m(1+x)}}\log^{k}\frac{e}{mx}\ dx\lesssim \int_{0}^{\frac{1}{2}}\frac{x^{n+\delta-t}}{\log^{k}\frac{e}{m}}\log^{k}\frac{e}{mx}\ dx
$$
$$
\lesssim \int_{0}^{\frac{1}{2}}x^{n+\delta-t}\ dx+\log^{-k}\frac{e}{m}\int_{0}^{\frac{1}{2}}x^{n+\delta-t}\log^{k}\frac{1}{x}\ dx\asymp 1.
$$

Otherwise, let $y=x/(x+1)$. We also have that
\begin{eqnarray*}
&\;&
\int_{\frac{1}{2}}^{\frac{1}{m}}\frac{x^{n+\delta-t}}{(1+x)^{r}\log^{k}\frac{e}{m(1+x)}}\log^{k}\frac{e}{mx}\ dx\\
&\;&
=\int_{\frac{1}{3}}^{\frac{1}{m+1}}y^{n+\delta-t}(1-y)^{r+t-\delta-n-2}\log^{k}\frac{e(1-y)}{my}\log^{-k}\frac{e(1-y)}{m}\ dy\\
&\;&
\lesssim \int_{\frac{1}{3}}^{1}y^{n+\delta-t}(1-y)^{r+t-\delta-n-2}\ dy\asymp 1.
\end{eqnarray*}
This shows that
$$
 J_{1}\lesssim  \frac{1}{|1-\langle w,\eta\rangle|^{r+t-\delta-n-1}} \ \mbox{when \ $t-\delta<n+1<r+t-\delta$ and $n>1$}.  \eqno{(3.5)}
$$

The proof of case $n=1$ is similar to the proof of Lemma 1 in [35].
\vskip2mm
It follows from (3.4) and (3.5) that the  $``\lesssim"$  part of (1) is  true. It follows from (3.3) and (3.4) that the  $``\lesssim"$  part of (2) is  true.   It follows from (3.2) and (3.4) that the  $``\lesssim"$  part of (3) is  true.
Next, we consider  $``\gtrsim"$  parts.

\vskip2mm
When $t-\delta<n+1$ and $r-\delta<n+1$, it follows from Proposition B(4) that
$$
J_{w,\eta}\gtrsim \int_{B_{n}}\frac{(1-|z|^{2})^{
\delta}}{|1-\langle z,w\rangle|^{ t}\ |1-\langle z,\eta\rangle|^{ r}}\
dv(z)\asymp \frac{1}{|1-\langle w,\eta\rangle|^{r+t-\delta-n-1}}.
$$
This shows that the  $``\gtrsim"$  part of (1) is  true.
\vskip2mm
When $t-\delta>n+1$, by  a change of variables $z=\varphi_{w}(u)$ and (1.1), we have that
\begin{eqnarray*}
&\;&
J_{w,\eta}= \frac{(1-|w|^{2})^{\delta+n+1-t}}{|1-\langle w,\eta\rangle|^{r}}\int_{B_{n}}\frac{ |1-\langle u,\varphi_{w}(\eta)\rangle|^{-r}(1-|u|^{2})^{
\delta}\log^{k}\frac{e|1-\langle u,w\rangle|^{2}}{(1-|u|^{2})(1-|w|^{2})}
dv(u)}{|1-\langle u,w\rangle|^{ 2\delta+2n+2-r-t}\log^{k}\frac{e|1-\langle u,\eta\rangle|}{|1-\langle w,\eta\rangle|\ |1-\langle u,\varphi_{w}(\eta)\rangle|}}\\
&\;&
\geq\frac{(1-|w|^{2})^{\delta+n+1-t}}{|1-\langle w,\eta\rangle|^{r}}\int_{|u|\leq\frac{1}{9}}\frac{ |1-\langle u,\varphi_{w}(\eta)\rangle|^{-r}(1-|u|^{2})^{
\delta}\log^{k}\frac{e|1-\langle u,w\rangle|^{2}}{(1-|u|^{2})(1-|w|^{2})}\
dv(u)}{|1-\langle u,w\rangle|^{ 2\delta+2n+2-r-t}\log^{k}\frac{e|1-\langle u,\eta\rangle|}{|1-\langle w,\eta\rangle|\ |1-\langle u,\varphi_{w}(\eta)\rangle|}}\\
&\;&
\asymp \frac{1}{(1-|w|^{2})^{t-\delta-n-1}|1-\langle w,\eta\rangle|^{r}}\log^{-k}\frac{e}{|1-\langle w,\eta\rangle|}\log^{k}\frac{e}{1-|w|^{2}}.
\end{eqnarray*}
This shows that the  $``\gtrsim"$  part of (3) is  true.
\vskip2mm
When $t-\delta=n+1>r-\delta$, by Lemma 1.8 in [1],  Lemma 2.2 and $(1-\rho|w|^{2})/2<1-\rho|w|\leq 1-\rho|w|^{2}$ for any $0\leq \rho<1$ and $w\in B_{n}$, we have that
$$
J_{w,\eta}\gtrsim \int_{0}^{1}\frac{(1-\rho^{2})^{\delta}}{(1-\rho^{2}|w|^{2})^{\delta+1}|1-\rho\langle w,\eta\rangle|^{r}}\log^{k}\frac{e}{1-\rho^{2}}\log^{-k}\frac{e}{|1-\rho\langle w,\eta\rangle|}\ d\rho
$$
$$ \gtrsim \int_{0}^{1}\frac{(1-\rho)^{\delta}}{(1-\rho|w|^{2})^{\delta+1}|1-\rho\langle w,\eta\rangle|^{r}}\log^{k}\frac{e}{1-\rho}\log^{-k}\frac{e}{|1-\rho\langle w,\eta\rangle|}\ d\rho.\eqno{(3.6)}
$$

When $|w|>1/\sqrt{2}$,  it follows from a change of variables \
$\displaystyle{x=\frac{|w|^{2}(1-\rho)}{1-\rho |w|^{2}}}$ \ that
\begin{eqnarray*}
&\;&
\int_{0}^{1}\frac{(1-\rho)^{\delta}}{(1-\rho|w|^{2})^{\delta+1}|1-\rho\langle w,\eta\rangle|^{r}}\log^{k}\frac{e}{1-\rho}\log^{-k}\frac{e}{|1-\rho\langle w,\eta\rangle|}\ d\rho \\
&\;&=
\frac{1}{|w|^{2(\delta+1)}}\int_{0}^{|w|^{2}}\frac{x^{\delta}}{1-x}\left|1-\langle
w,\eta\rangle+\frac{\langle
w,\eta\rangle(1-|w|^{2})x}{|w|^{2}(1-x)}\right|^{-r}\log^{k}\frac{e|w|^{2}(1-x)}{(1-|w|^{2})x}\\
&\;&
\times \ \log^{-k}\frac{e}{\left|1-\langle
w,\eta\rangle+\frac{\langle
w,\eta\rangle(1-|w|^{2})x}{|w|^{2}(1-x)}\right|}\ dx
\\
&\;& \gtrsim\int_{\frac{1}{2}}^{|w|^{2}}\frac{1}{1-x}\left(m+\frac{1-|w|^{2}}{1-x}\right)^{-r}\log^{k}\frac{e(1-x)}{1-|w|^{2}}\log^{-k}\frac{e}{m+\frac{1-|w|^{2}}{1-x}}\
dx\\
&\;&
 \gtrsim \log^{-k}\frac{e}{m}\int_{\frac{1}{2}}^{|w|^{2}}\frac{1}{1-x}\left(m+\frac{1-|w|^{2}}{1-x}\right)^{-r}\log^{k}\frac{e(1-x)}{1-|w|^{2}}\
dx.
\end{eqnarray*}

It follows from (3.12) and (3.13) in [12] that
$$
\int_{0}^{1}\frac{(1-\rho)^{\delta}}{(1-\rho|w|^{2})^{\delta+1}|1-\rho\langle w,\eta\rangle|^{r}}\log^{k}\frac{e}{1-\rho}\log^{-k}\frac{e}{|1-\rho\langle w,\eta\rangle|}\ d\rho
$$
$$
\gtrsim \frac{1}{m^{r}}\log^{-k}\frac{e}{m}\log\frac{em}{1-|w|^{2}}\log^{k}\frac{e}{1-|w|^{2}}. \ \ \ \ \ \ \ \ \ \ \ \ \ \ \ \ \ \ \ \ \ \ \ \eqno{(3.7)}
$$

When $t-\delta=n+1>r-\delta$, it follows from  (3.6) and (3.7) that
$$
 J_{w,\eta}\gtrsim \frac{1}{|1-\langle w,\eta\rangle|^{r}}\log^{-k}\frac{e}{|1-\langle w,\eta\rangle|}\log^{k}\frac{e}{1-|w|^{2}}\log\frac{e}{|1-\langle w,\varphi_{w}(\eta)\rangle|}.
 $$
This shows that the  $``\gtrsim"$  part of (2) is  true.
This proof is complete. \ \ $\Box$

\vskip2mm
{\bf Proposition 3.3} \ Let $\mu_{1}$, $\mu_{2}$ be two normal functions on $[0,1)$, $h\in H(B_{n})$. If
$$
\sup_{z\in B_{n}}\frac{\mu_{1}(|z|)|h(z)|}{\mu_{2}(|z|)}=M<\infty, \ \ \mbox{then}
$$
$$
\sup_{z\in B_{n}}\frac{(1-|z|^{2})\mu_{1}(|z|)|\nabla h(z)|}{\mu_{2}(|z|)}\lesssim M.
$$

{\bf Proof} \ Let $g$ be this function in Lemma 2.5, matching $\mu_ {2}$. For any $w\in B_{n}$, let $
F_{w}(z)=h(z)g(\langle z,w\rangle)$. It follows from Lemma 2.5 that
$$
\mu_{1}(|z|)|F_{w}(z)|=\frac{\mu_{1}(|z|)|h(z)|}{\mu_{2}(|z|)}\mu_{2}(|z|)|g(\langle z,w\rangle)|\leq M\mu_{2}(|\langle z,w\rangle|)|g(\langle z,w\rangle)|\leq M_{1}M
$$
for all $z\in B_{n}$. Moreover, we have that
$|F_{w}(z)|\leq M_{1}M/\mu_{1}(0)(1-|z|^{2})^{b_{1}} $, where $b_{1}$  the larger parameter in the definition of $\mu_{1}$. When $\alpha>b_{1}-1$, it is clear that $F_{w}\in A_{\alpha}^{1}(B_{n})$. It follows from Lemma 2.6 that
$$
F_{w}(z)=\int_{B_{n}}\frac{F_{w}(u)\ dv_{\alpha}(u)}{(1-\langle z,u\rangle)^{n+1+\alpha}} \ \ (z\in B_{n}).
$$

Let  $a_{1}$ be the smaller parameter in the definition of $\mu_{1}$. By Lemma 2.3 and Proposition A, we have that
\begin{eqnarray*}
 &\;&
\mu_{1}(|z|)|(\nabla F_{w})(z)|\leq \int_{B_{n}}\frac{(n+1+\alpha)M_{1}Mc_{\alpha}\mu_{1}(|z|)(1-|u|^{2})^{\alpha} dv(u)}{\mu_{1}(|u|)|1-\langle z,u\rangle|^{n+2+\alpha}}\\&\;&
\lesssim  \int_{B_{n}}\frac{M(1-|z|^{2})^{a_{1}}(1-|u|^{2})^{\alpha-a_{1}} dv(u)}{|1-\langle z,u\rangle|^{n+2+\alpha}}+\int_{B_{n}}\frac{M(1-|z|^{2})^{b_{1}}(1-|u|^{2})^{\alpha-b_{1}} dv(u)}{|1-\langle z,u\rangle|^{n+2+\alpha}}\\
&\;&
\asymp \frac{M}{1-|z|^{2}}.
\end{eqnarray*}

In particular, if we take $z=w$,  then it is clear  that
$$
(1-|w|^{2})\mu_{1}(|w|)|g(|w|^{2})\nabla h(w)+h(w)g'(|w|^{2})\overline{w}|\lesssim M.
$$

Suppose that $b_{2}$ is the larger parameter in the definition of $\mu_{2}$. By Lemma 2.5 and the definition of normal function,  it is clear that
$$
\frac{(1-|w|^{2})\mu_{1}(|w|)|\nabla h(w)|}{\mu_{2}(|w|)}\leq\frac{2^{b_{2}}}{M_{0}}  (1-|w|^{2})\mu_{1}(|w|)|\nabla h(w)|g(|w|^{2})
$$
$$
\lesssim M+\frac{\mu_{1}(|w|)|h(w)|}{\mu_{2}(|w|)}(1-|w|^{2})\mu_{2}(|w|)g'(|w|^{2})\leq M(1+M_{2}).
$$

This proof is complete.  \ \ $\Box$
\vskip2mm
{\bf Theorem 3.4} \ Let $\mu$ and $\nu$ be two normal functions on $[0,1)$. Suppose that $\psi\in H(B_{n})$ and  $0\leq s\leq n$. Then there are the following results:
\vskip3mm
(1) \ If $\int_{0}^{1}\mu^{-1}(\rho)(1-\rho^{2})^{\frac{s-n}{p}}\ d\rho<\infty$, then $\psi$
is a pointwise  multiplier from $F(p,\mu,s)$ to $ \mathcal{B_{\nu}}(B_{n})$ if and only if $\psi\in \mathcal{B_{\nu}}(B_{n})$ and
$$
\sup_{z\in B_{n}}\frac{\nu(|z|)|\psi(z)|}{\mu(|z|)(1-|z|^{2})^{\frac{n-s}{p}}}<\infty. \eqno{(3.8)}
$$

(2) \ If $\int_{0}^{1}\mu^{-1}(\rho)(1-\rho^{2})^{\frac{s-n}{p}}\ d\rho=\infty$ and $1+\int_{0}^{r}\mu^{-1}(\rho)(1-\rho^{2})^{\frac{s-n}{p}}\ d\rho\asymp \mu^{-1}(r)(1-r^{2})^{1+\frac{s-n}{p}}$,
then $\psi$ is a  pointwise multiplier from $F(p,\mu,s)$ to $ \mathcal{B_{\nu}}(B_{n})$ if and only if
(3.8) holds.

\vskip2mm
{\bf Proof} \ Suppose that $\psi\in \mathcal{B_{\nu}}(B_{n})$ and (3.8) holds.
For any  $f\in F(p,\mu,s)$ and $z\in B_{n}$ with $|z|>1/2$, it follows from Lemma 2.4 and Proposition 3.3 that
\begin{eqnarray*}
 &\;&
 \nu(|z|)|\nabla (\psi f)(z)|\leq  \nu(|z|)|\psi(z)||\nabla f(z)|+\nu(|z|)|\nabla \psi(z)|| f(z)|
 \\
 &\;&
 \lesssim \frac{ \nu(|z|)|\psi(z)|}{(1-|z|^{2})^{\frac{n-s}{p}}\mu(|z|)}\ ||f||_{p,\mu,s}\\
 &\;&+ \ \nu(|z|)|\nabla \psi(z)|\left\{ 1+\int_{0}^{|z|}\frac{d\rho}{\mu(\rho)(1-\rho^{2})^{\frac{n-s}{p}}}\right\}||f||_{p,\mu,s}\\
 &\;&
 \lesssim ||f||_{p,\mu,s} \ \ \Rightarrow \ \ ||\psi f||_{p,\mu,s}<\infty.
 \end{eqnarray*}
This means that $\psi$ is a pointwise multiplier from $F(p,\mu,s)$ to $ \mathcal{B_{\nu}}(B_{n})$.

\vskip2mm
Conversely, let $\psi$ be a pointwise multiplier from $F(p,\mu,s)$ to $ \mathcal{B_{\nu}}(B_{n})$. It follows from the closed graph theorem that
$M_{\psi}: f\rightarrow \psi f$ is a bounded operator from $F(p,\mu,s)$ to $ \mathcal{B_{\nu}}(B_{n})$. First, it is clear that $\psi\in \mathcal{B_{\nu}}(B_{n})$ by taking $f(z)=1$. Suppose that $b$ is the larger parameter in the definition of $\mu$. For any $w\in B_{n}$ with $|w|>1/2$, we take that
$$
f_{w}(z)=\frac{(1-|w|^{2})^{1+b+\frac{s}{p}}}{\mu(|w|)}\left\{\frac{1-|w|^{2}}{(1-\langle z,w\rangle)^{b+1+\frac{n}{p}}}-\frac{1}{(1-\langle z,w\rangle)^{b+\frac{n}{p}}}\right\} \ \ \ (z\in B_{n}).
$$
Then $f_{w}(w)=0$ and it follows from Lemma 2.3 that
\begin{eqnarray*}
 &\;&
 \int_{B_{n}}\frac{(1-|u|^{2})^{s}|\nabla f_{w}(z)|^{p}\ \mu^{p}(|z|)}{|1-\langle z,u\rangle|^{2s}(1-|z|^{2})}\ dv(z)\\
 &\;&
   \lesssim\int_{B_{n}}\frac{(1-|u|^{2})^{s}(1-|w|^{2})^{p+pb+s}\ \mu^{p}(|z|)}{\mu^{p}(|w|)|1-\langle z,u\rangle|^{2s}|1-\langle z,w\rangle|^{p+pb+n}(1-|z|^{2})}\ dv(z)
  \\
 &\;&
 \lesssim \int_{B_{n}}\frac{(1-|u|^{2})^{s}(1-|w|^{2})^{p+s}(1-|z|^{2})^{pb-1}}{|1-\langle z,u\rangle|^{2s}|1-\langle z,w\rangle|^{p+pb+n}}\ dv(z)\\&\;&
 + \ \int_{B_{n}}\frac{(1-|u|^{2})^{s}(1-|w|^{2})^{p+s+pb-pa}(1-|z|^{2})^{pa-1}}{|1-\langle z,u\rangle|^{2s}|1-\langle z,w\rangle|^{p+pb+n}}\ dv(z)=J_{1}+J_{2}.
 \end{eqnarray*}

When $0\leq 2s<pb+n$,  it follows from Proposition B(6) or Proposition A  that
\begin{eqnarray*}
 &\;&
 J_{1}\asymp \frac{(1-|u|^{2})^{s}(1-|w|^{2})^{p+s}}{(1-|w|^{2})^{p}|1-\langle u,w\rangle|^{2s}}\lesssim 1.
 \end{eqnarray*}

When  $2s>pb+n$,  it follows from Proposition B(7) and $0\leq s\leq n$ that
\begin{eqnarray*}
 &\;&
  J_{1}\asymp\frac{(1-|u|^{2})^{s}(1-|w|^{2})^{p+s}}{(1-|w|^{2})^{p}|1-\langle u,w\rangle|^{2s}}+\frac{(1-|u|^{2})^{s}(1-|w|^{2})^{p+s}}{(1-|u|^{2})^{2s-pb-n}|1-\langle u,w\rangle|^{p+pb+n}}\lesssim 1.
 \end{eqnarray*}

 When $2s=pb+n$, we choose  $0<\varepsilon<s$. By Proposition B(8) and (2.1), we have
\begin{eqnarray*}
  &\;& J_{1}\asymp\frac{(1-|u|^{2})^{s}(1-|w|^{2})^{p+s}}{(1-|w|^{2})^{p}|1-\langle u,w\rangle|^{2s}} + \ \frac{(1-|u|^{2})^{s}(1-|w|^{2})^{p+s}}{|1-\langle u,w\rangle|^{p+pb+n}}\log\frac{e}{|1-\langle \varphi_{u}(w), u\rangle|}
  \end{eqnarray*}
  \begin{eqnarray*}
  &\;&
  \lesssim \frac{(1-|u|^{2})^{s}(1-|w|^{2})^{p+s}}{(1-|w|^{2})^{p}|1-\langle u,w\rangle|^{2s}}+ \frac{(1-|u|^{2})^{s-\varepsilon}(1-|w|^{2})^{p+s}}{|1-\langle u,w\rangle|^{p+pb+n-\varepsilon}}\lesssim 1.
 \end{eqnarray*}

Similarly, we may prove that $J_{2}\lesssim 1$.
 This shows that $||f_{w}||_{p,\mu,s}\lesssim 1$.
\vskip2mm
If $\psi$ is a  pointwise multiplier from $F(p,\mu,s)$ to $ \mathcal{B_{\nu}}(B_{n})$, then we have that
 $$
 ||M_{\psi}||\gtrsim ||M_{\psi}||.||f_{w}||_{p,\mu,s}\geq \nu(|w|)|\nabla(\psi f_{w})(w)|\asymp\frac{\nu(|w|)|\psi(w)|}{(1-|w|^{2})^{\frac{n-s}{p}}\mu(|w|)}.
$$
This shows that (3.8) holds. \ This proof is complete. \ \ $\Box$
\vskip2mm

{\bf Corollary 3.5} \ Let $\mu$ and $\nu$ be two normal functions on $[0,1)$. Suppose that $\psi\in H(B_{n})$, $s>n$ and $\gamma(\rho)=(1-\rho^{2})^{\frac{n-s}{p}}\mu(\rho)$ is a normal function on $[0,1)$. Then $\psi$ is a  pointwise multiplier from $F(p,\mu,s)$ to $ \mathcal{B_{\nu}}(B_{n})$ if and only if
(3.8) holds and
$$
\sup_{z\in B_{n}}\nu(|z|)|\nabla \psi(z)|\int_{0}^{|z|}\frac{d\rho}{\mu(\rho)(1-\rho^{2})^{\frac{n-s}{p}}}<\infty.
$$

{\bf Proof} \ It follows from Lemma 2.4 that $F(p,\mu,s)=\mathcal{B_{\gamma}}(B_{n})$ when $s>n$. By Theorem A in [6], we may obtain this result.

\vskip2mm
In Theorem 3.4, if $\int_{0}^{1}\mu^{-1}(\rho)(1-\rho^{2})^{\frac{s-n}{p}} d\rho=\infty$ and $ 1+\int_{0}^{r}\mu^{-1}(\rho)(1-\rho^{2})^{\frac{s-n}{p}} d\rho$  is not equivalent to $ \mu^{-1}(r)(1-r^{2})^{1+\frac{s-n}{p}}$ for  $0\leq r<1$, then we do not give the the necessary and sufficient conditions.
However, there are the following results for the normal weight $\mu(r)=(1-r^{2})^{\alpha}\log^{\beta}\frac{e}{1-r^{2}}$ \ ($\alpha>0$, $-\infty<\beta<\infty$).
\vskip2mm
{\bf Corollary 3.6} \ Let $\nu$ be a normal function on $[0,1)$. Suppose that $\mu(r)=(1-r^{2})^{\alpha}\log^{\beta}\frac{e}{1-r^{2}}$ \ ($\alpha>0$, $-\infty<\beta<\infty$), $\psi\in H(B_{n})$ and $0\leq s\leq n$. Then there are the following results:
\vskip2mm
(1) \ If \ $p\alpha+n-p<s$, or \ $p\alpha+n-p=s>0$ and $\beta>1$, then $\psi$ is a  pointwise multiplier from $F(p,\mu,s)$ to $ \mathcal{B_{\nu}}(B_{n})$
if and only if $\psi\in \mathcal{B_{\nu}}(B_{n})$ and (3.8) holds.
\vskip2mm
(2) \ If \ $p\alpha+n-p>s$, then $\psi$ is a  pointwise multiplier from $F(p,\mu,s)$ to $ \mathcal{B_{\nu}}(B_{n})$
if and only if (3.8) holds.
\vskip2mm
(3) \ If \ $p\alpha+n-p=s>0$ and $\beta=1$, then $\psi$ is a  pointwise multiplier from $F(p,\mu,s)$ to $ \mathcal{B_{\nu}}(B_{n})$
if and only if (3.8) holds and
$$
\sup_{z\in B_{n}}\nu(|z|)|\nabla \psi(z)|\log\log\frac{e}{1-|z|^{2}}<\infty.\eqno{(3.9)}
$$

(4) \ If \ $p\alpha+n-p=s>0$ and  $\beta<1$, then $\psi$ is a  pointwise multiplier from $F(p,\mu,s)$ to $ \mathcal{B_{\nu}}(B_{n})$
if and only if (3.8) holds and
$$
\sup_{z\in B_{n}}\nu(|z|)|\nabla \psi(z)|\log^{1-\beta}\frac{e}{1-|z|^{2}}<\infty.\eqno{(3.10)}
$$

 {\bf Proof} \ The necessary and sufficient conditions of (1)-(2), and the sufficient conditions of (3)-(4) come from Theorem 3.4.
\vskip2mm
For any $w\in B_{n}$ with $|w|>1/2$, we take the test function
$$
G_{w}(z)= \left\{\int_{0}^{\langle z,w\rangle}\frac{1}{1-t}\log^{-\beta}\frac{e}{1-t} \ dt\right\}^{2}/\int_{0}^{|w|^{2}}\frac{1}{1-t}\log^{-\beta}\frac{e}{1-t} \ dt
$$$$- \ 2\int_{0}^{\langle z,w\rangle}\frac{1}{1-t}\log^{-\beta}\frac{e}{1-t} \ dt \ \ \ (z\in B_{n}).
$$

We know that $p+s-n-1=p\alpha-1>-1$, \ $2s+p-(p+s-n-1)=s+n+1>n+1$ and
$p-(p+s-n-1)=-s+n+1<n+1$. When  $\beta>0$ and $0\leq s<p$, for any  $u\in B_{n}$, it follows from Proposition 3.2(1) that
\begin{eqnarray*}
 &\;&
\int_{B_{n}}\frac{(1-|u|^{2})^{s}|\nabla G_{w}(z)|^{p}\ \mu^{p}(|z|)\ dv(z)}{(1-|z|^{2})|1-\langle z,u\rangle|^{2s}}
 \\
&\;&\lesssim \int_{B_{n}}\frac{(1-|u|^{2})^{s}(1-|z|^{2})^{p+s-n-1}\log^{p\beta}\frac{e}{1-|z|^{2}}\ dv(z)}{|1-\langle z,u\rangle|^{2s}|1-\langle z,w\rangle|^{p}\log^{p\beta}\frac{e}{|1-\langle z,w\rangle|}}\\
&\;&
\asymp \frac{(1-|u|^{2})^{s}}{|1-\langle u,w\rangle|^{2s+p-(p+s-n-1)-n-1}}=\frac{(1-|u|^{2})^{s}}{|1-\langle u,w\rangle|^{s}}
\lesssim 1.
\end{eqnarray*}

When  $\beta>0$ and $s=p$, we choose $0<\varepsilon<p$.  For any  $u\in B_{n}$,  it follows from  Proposition 3.2(2) and (2.1) that
\begin{eqnarray*}
 &\;&
\int_{B_{n}}\frac{(1-|u|^{2})^{s}|\nabla G_{w}(z)|^{p}\ \mu^{p}(|z|)\ dv(z)}{(1-|z|^{2})|1-\langle z,u\rangle|^{2s}}\\
&\;&
 \lesssim \frac{(1-|u|^{2})^{s}}{|1-\langle u,w\rangle|^{s}}\log^{-p\beta}\frac{e}{|1-\langle u,w\rangle|}\log^{p\beta}\frac{e}{1-|u|^{2}}\log\frac{e}{|1-\langle u,\varphi_{u}(w)\rangle|}\\
 &\;&
 =\frac{(1-|u|^{2})^{s-\varepsilon}}{|1-\langle u,w\rangle|^{s-\varepsilon}}|1-\langle u,\varphi_{u}(w)\rangle|^{\frac{\varepsilon}{2}}\log\frac{e}{|1-\langle u,\varphi_{u}(w)\rangle|}\\
 &\;&
 \times \ (1-|u|^{2})^{\frac{\varepsilon}{2}}\log^{p\beta}\frac{e}{1-|u|^{2}}|1-\langle u,w\rangle|^{-\frac{\varepsilon}{2}}\log^{-p\beta}\frac{e}{|1-\langle u,w\rangle|}\lesssim 1.
\end{eqnarray*}

When  $\beta>0$ and $s>p$,  it follows from Proposition 3.2(3) and (2.1) that
\begin{eqnarray*}
 &\;&
\int_{B_{n}}\frac{(1-|u|^{2})^{s}|\nabla G_{w}(z)|^{p}\ \mu^{p}(|z|)\ dv(z)}{(1-|z|^{2})|1-\langle z,u\rangle|^{2s}}\\
&\;&
 \lesssim \frac{(1-|u|^{2})^{s}}{(1-|u|)^{s-p}|1-\langle u,w\rangle|^{p}}\log^{-p\beta}\frac{e}{|1-\langle u,w\rangle|}\log^{p\beta}\frac{e}{1-|u|^{2}}\lesssim 1.
\end{eqnarray*}

When  $\beta\leq 0$, it follows from Proposition B(6)-(8) that
\begin{eqnarray*}
 &\;&
\int_{B_{n}}\frac{(1-|u|^{2})^{s}|\nabla G_{w}(z)|^{p}\ \mu^{p}(|z|)\ dv(z)}{(1-|z|^{2})|1-\langle z,u\rangle|^{2s}}\\
&\;&\lesssim \int_{B_{n}}\frac{(1-|u|^{2})^{s}(1-|z|^{2})^{p+s-n-1}\ dv(z)}{|1-\langle z,u\rangle|^{2s}|1-\langle z,w\rangle|^{p}}\lesssim 1.
\end{eqnarray*}

This shows that $||G_{w}||_{p,\mu,s}\lesssim 1$.
If $\psi$ is a pointwise multiplier from $F(p,\mu,s)$ to $ \mathcal{B_{\nu}}(B_{n})$, then we may obtain that
$$
\nu(|w|)|\nabla\psi(w)|\int_{0}^{|w|^{2}}\frac{1}{1-t}\log^{-\beta}\frac{e}{1-t} \ dt\lesssim ||M_{\psi}||,
$$

According to the arbitrariness of $w$,  we can get (3.9) and (3.10) by different cases.

\vskip2mm  This proof is complete. \ \ $\Box$
\vskip2mm
In particular, when $\nu(r)=(1-r^{2})^{\frac{n-s}{p}}\mu(r)$ in Corollary 3.6, it follows from Lemma 2.4 that $F(p,\mu,s)\subseteq  \mathcal{B_{\nu}}(B_{n})$. Moreover, there are the following results:
\vskip2mm
{\bf Corollary 3.7} \ Let \ $\mu(r)=(1-r^{2})^{\alpha}\log^{\beta}\frac{e}{1-r^{2}}$.  If \ $\nu(r)=(1-r^{2})^{\alpha+\frac{n-s}{p}}\log^{\beta}\frac{e}{1-r^{2}}$ \ ($\alpha>0$, $-\infty<\beta<\infty$) is a normal function on $[0,1)$ and $\psi\in H(B_{n})$, then   $\psi$ is a pointwise multiplier from $F(p,\mu,s)$ to $ \mathcal{B_{\nu}}(B_{n})$
if and only if
\vskip2mm
(1) \ $\psi\in \mathcal{B_{\nu}}(B_{n})$ \ when \ $p\alpha+n-p<s$,  or \ $p\alpha+n-p=s>0$ and $\beta>1$.
\vskip2mm
(2) \ $\psi\in H^{\infty}(B_{n})$ \ when \ $p\alpha+n-p>s$.
\vskip2mm
(3) \  $\psi\in H^{\infty}(B_{n})$ \ and
$$\displaystyle{\sup_{z\in B_{n}}(1-|z|^{2})|\nabla \psi(z)|\left(\log\frac{e}{1-|z|^{2}}\right)\log\log\frac{e}{1-|z|^{2}}<\infty}\eqno{(3.11)}$$
when \ $p\alpha+n-p=s>0$ and $\beta=1$.
\vskip2mm
(4) \  $\psi\in H^{\infty}(B_{n})$ \ and
$$\displaystyle{\sup_{z\in B_{n}}(1-|z|^{2})|\nabla \psi(z)|\log\frac{e}{1-|z|^{2}}<\infty}\eqno{(3.12)}$$
when \  $p\alpha+n-p=s>0$ and $\beta<1$.
\vskip2mm
{\bf Proof} \ It follows from Lemma 2.4 that $\mathcal{B_{\nu}}(B_{n})\subseteq H^{\infty}(B_{n})$ when $\int_{0}^{1}\mu^{-1}(\rho)(1-\rho^{2})^{\frac{s-n}{p}}<\infty$. By  Corollary 3.5 and Corollary 3.6, we may obtain these results.
\vskip2mm
{\bf Note 2} \ In Corollary 3.7, no condition can be removed.  For example, if  $\psi(z)=e^{\frac{z_{1}+1}{z_{1}-1}}$, then
 $\psi\in  H^{\infty}(B_{n})$.   But
\begin{eqnarray*}
 &\;&
\sup_{z\in B_{n}}(1-|z|^{2})|\nabla \psi(z)|\log\frac{e}{1-|z|^{2}}\geq \sup_{z_{1}\in D}(1-|z_{1}|^{2})|\nabla \psi(z_{1})|\log\frac{e}{1-|z_{1}|^{2}}\\
&\;&
=\sup_{re^{i\theta}\in D}\frac{2(1-r^{2})}{1+r^{2}-2r\cos\theta}e^{\frac{r^{2}-1}{1+r^{2}-2r\cos\theta}}\log\frac{e}{1-r^{2}}\geq\sup_{0<\theta<\frac{\pi}{2}}\frac{2}{e}\log\frac{e}{\sin^{2}\theta}=\infty.
\end{eqnarray*}
This means that $\psi$ does not meet (3.11) and (3.12). If $\psi(z)=\log\log\frac{e^{2}}{1-z_{1}}$, then it meets (3.12). But $\psi$ does not belong to $H^{\infty}(B_{n})$. If $\psi(z)=\log\log\log\frac{e^{4}}{1-z_{1}},$ then it meets (3.11). But $\psi$ does not belong to $H^{\infty}(B_{n})$.

\vskip2mm
{\bf Statement:}  No potential conflict of interest.

\end{document}